\documentclass[12pt]{article}

\usepackage{latexsym,amsmath,amssymb}

\usepackage{algorithmic}
\usepackage{algorithm}
\usepackage{graphicx}
\usepackage{eepic}
\input{epsf.sty}

\title{REGULAR CLOSED SETS OF PERMUTATIONS}

\author{M.H. Albert, M.D. Atkinson\\
Department of Computer Science\\
University of Otago, New Zealand\\and \\
N. Ru\v{s}kuc\\
School of Mathematics and Statistics\\
University of St Andrews, UK}

\renewcommand{\to}{\rightarrow}

\renewcommand{\int}{\cap}
\newcommand{\ci}{\subseteq}

\newtheorem{lemma}{Lemma}
\newtheorem{theorem}{Theorem}

\newtheorem{proposition}{Proposition}
\newtheorem{corollary}{Corollary}

\newtheorem{exx}{Example}
\newenvironment{example}{\begin{exx}\rm}{\end{exx}}
\newtheorem{remm}{Remark}

\newtheorem{noo}{Note}
\newenvironment{note}{\begin{noo}\rm}{\end{noo}}
\newcommand{\qed}{\hfill \rule{1ex}{1ex}} 

\newenvironment{proof}{{\bf Proof}: }{\qed}

\renewcommand{\L}{{\cal L}}

\newcommand{\R}{{\cal R}}
\newcommand{\B}{{\cal B}}
\newcommand{\W}{{\cal W}}
\newcommand{\X}{{\cal X}}

\begin{document}

\maketitle

\begin{abstract}Machines whose main purpose is to permute and sort
data are studied.  The sets of permutations that can arise are
analysed by means of finite automata and avoided pattern techniques.
Conditions are given for these sets being enumerated by rational
generating functions.
\end{abstract}


\section{Introduction}

From the earliest days of Computer Science abstract machines have
been used to model computations and categorise them according to the
different resources they require.  In this paper we consider a new
type of machine that is suited to modelling computations whose sole or
main effect is to permute data.  Unlike most classical machines these
new machines have an infinite input alphabet whose symbols form the
data that is to be permuted.  Despite this we shall show how the
theory of finite automata can be deployed in their analysis.

A \emph{permuting machine} is a non-deterministic machine with the
following properties:

\begin{enumerate}
    \item it transforms an input stream of distinct tokens
    into an output stream that is
    a permutation of the input stream,
    \item it is oblivious to the values of the input stream tokens,
    \item it has a hereditary property: if an input stream $\sigma$ can be
    transformed into an output stream $\tau$ and $\sigma'$ is a
    subsequence of $\sigma$ whose symbols transform into the
    subsequence $\tau'$ of $\tau$, then it is possible for $\sigma'$
    (if presented as an input to the machine in its own right) to be
    transformed into an output stream $\tau'$.
\end{enumerate}

\noindent\textbf{Examples}

\begin{enumerate}
\item A \emph{riffle shuffler} divides the input stream into two
segments and then interleaves them
in any way to form the output stream.
\item A \emph{stack} receives members of the input stream and
outputs them under a last-in-first-out discipline.
\item A \emph{transportation network} \cite{AtkLivTul} is any finite
directed graph
with a node to represent the input stream and a node to represent the
output stream.  The other nodes can each hold one of the input objects
and the objects are moved around the graph until they emerge at the
output node.
\end{enumerate}

The oblivious property of permuting machines allows us to name the
input tokens $1,2,\ldots,n$ (in that order) in which case the output
will be some permutation of  $1,2,\ldots,n$.  In this way we can
consider a permuting machine to be a generator of permutations (usually,
because of the non-determinism, generating many of each length).  There is
another point of view which is sometimes more useful where we
consider the input stream to be some permutation of $1,2,\ldots,n$ and
ask whether the machine is capable of sorting the tokens (so that they
appear in the output stream in the order  $1,2,\ldots,n$).  These two
viewpoints are equivalent since a machine can generate a particular
permutation $\sigma$ if and only if it can sort the permutation
$\sigma^{-1}$.

However, it is the hereditary property which allows non-trivial
properties of permuting machines to be found because of a connection
with the combinatorial theory of involvement and closed sets of
permutations.  Formally, a permutation
$\pi$ is said to be \emph{involved} in
another permutation $\sigma$ (denoted as $\pi\preceq\sigma$)
if $\pi$ is order isomorphic to a
subsequence of $\sigma$.  For example $231$ is involved in $31542$
because of the subsequence $352$ (or the subsequence $342$).  We also
say that $\sigma$ \emph{avoids} $\pi$ if $\pi$ is not involved in
$\sigma$.

Permutation involvement has been an active area of combinatorics for
over 10 years although it surfaced long before that in data
structuring questions on stacks, queues and deques (see \cite{Knu,
Pra, Tar}).  Involvement is a partial order on the set of all
permutations and is conveniently studied by means of order ideals
called closed sets.  A \emph{closed} set $\mathcal{X}$
of permutations is one  with the property that $\sigma\in\mathcal{X}$
and $\pi\preceq\sigma$ imply $\pi\in\mathcal{X}$.  The connection
between permuting machines and closed sets is via the following
result which follows from the definitions.

\begin{proposition}
    The set of permutations that a permuting machine can generate,
    and the set that it can sort, are both closed.
\end{proposition}

In classical automata theory machines are associated with the
languages they recognise.  The above proposition suggests that the
appropriate associated language of a permuting machine is the closed
set of permutations that it can generate.  We will study
various permuting machines and their associated closed sets, and will
show the utility of the
permuting machine paradigm as a tool for advancing the theory of
permutation involvement.  Before giving further details of our
results we recall some key concepts about permutation involvement.

A closed set $\mathcal{X}$ is, by definition, closed
``downwards''.  But that is equivalent to its complement
$\mathcal{X}^{C}$ being closed ``upwards''
($\sigma\in\mathcal{X}^{C}$ and $\sigma\preceq\tau$ imply
$\tau\in\mathcal{X}^{C}$).  Obviously, $\mathcal{X}^{C}$ is
determined by its set of minimal permutations which we denote by
$B(\mathcal{X})$ and call the \emph{basis} of $\mathcal{X}$.  Clearly
\[\mathcal{X}=\{\sigma\mid\sigma\not\in\mathcal{X}^{C}\}=
\{\sigma\mid\pi\not\preceq\sigma\mbox{ for all }\pi\in
B(\mathcal{X})\}\]
is determined by its basis.  By definition, $B(\mathcal{X})$ is an
antichain in the involvement order and conversely every antichain has
the form $B(\mathcal{X})$ for some closed set $\mathcal{X}$.  The
bases of the closed sets of permutations generated by the
machines in examples $1$ and $2$ above are $\{321,2413,2143\}$ and
$\{312\}$ respectively.
The closed sets that arise in practice are generally
infinite so it is clearly significant to know when a finite
description is available by means of the basis.  Indeed, many
combinatorial enumeration investigations begin from some particular
finite basis and study properties of the closed set that it defines
(\cite{Bon,Wes}).  We let $\mathcal{A}(B)$ denote the closed set
whose basis is the antichain $B$; in other words
\[\mathcal{A}(B)=\{\sigma\mid\beta\not\preceq\sigma\mbox{ for all
}\beta\in B\}\]

Given a closed set $\mathcal{X}$ (or a permuting machine that defines
it) we would like to be able to solve

\begin{itemize}
    \item The decision problem: given a permutation $\sigma$ decide
    whether $\sigma\in\mathcal{X}$ (in linear time if possible),
    \item The enumeration problem: determine,
    for each length $n$, the number of permutations in $\mathcal{X}$,
    \item The basis problem: find the basis of $\mathcal{X}$, or at
    least determine whether the basis is finite or infinite.
\end{itemize}

In this paper we shall show how to exploit the classical
theory of finite automata
to make progress on these problems.
To do this we have to overcome the difficulty that
this theory deals with strings over a finite alphabet,
whereas the strings of $\mathcal{X}$ are written in the infinite
alphabet $1,2,\ldots$.  Therefore we shall
look for encodings of the permutations in
$\mathcal{X}$ as strings over a finite alphabet (normally
$[k]=\{1,2,\ldots,k\}$) and hope to prove that
the language of such encodings is regular (or to find conditions under
which this is so).  Once we have proved the regularity of such a
language we can appeal to two well-known facts: that regular languages
have linear time recognisers, and that the generating
function (the formal power series whose
coefficients give the number of sequences of each length) is a rational
function.

Of course this approach cannot be expected to succeed in
all cases if for no other reason than that closed sets do not always
have rational generating functions.  Nevertheless, in Sections
\ref{BoundedSection} and \ref{MonotoneSection}, we shall give two wide
classes of closed
sets (and permuting machines) which show that the approach can have
significant successes.  In particular we produce infinite families of
closed sets all of whose finitely based closed subsets have rational
generating functions.  Our results therefore link to the many recent
papers where particular closed sets have been enumerated (for
example, \cite{Bar,Cho,Man}).  In the final section we indicate how we hope
our approach may be extended.

We conclude this section by recalling some basic facts about
transducers.

For our purposes a transducer is essentially a (non-deterministic) finite
    automaton with output symbols (from an alphabet $\Gamma$)
    as well as input symbols (from an alphabet $\Delta$).  We
    allow $\epsilon$ inputs as well as $\epsilon$ outputs.
    A transducer defines a relation between $\Delta^{\ast}$ and
    $\Gamma^{\ast}$ in a natural way.  That is to say, for every path
    in the transducer from the starting state to one of the final
    states, let the sequence of input labels be $\alpha$ and the
    sequence of output labels $\beta$ (all $\epsilon$'s being omitted
    of course); then $(\alpha,\beta)$ is a related pair.

    In any transducer we can interchange the input and output symbols
    on each transition to obtain another transducer.  Therefore

    \begin{lemma} \label{transpose}
	If $\R$ is a transducer relation so also is the
        transpose relation $\R^{t}$.
    \end{lemma}

    Let $\mathcal{L}\subseteq \Delta^{\ast}$ and define
    \[\mathcal{L}\R=\{\beta\in\Gamma^{\ast}\mid \mbox{ there exists
    }\alpha\in \mathcal{L}\mbox{ with }(\alpha,\beta)\in \R\}\]

The main result we need from the theory of transducers appears as
exercise 11.9 in \cite{HopUll}.  For completeness, and to establish
notation, we include the proof.

\begin{proposition} \label{technical}
	If $\R$ is a transducer relation and $\mathcal{L}$ is a regular subset
	of $\Delta^{\ast}$ then $\mathcal{L}\R$ is regular.
    \end{proposition}
    \begin{proof}
	Let $P$ be the set of states of the transducer, $\mu$ the transition
	function (mapping $P\times(\Delta\cup\{\epsilon\})$ into subsets of
	$P\times(\Gamma\cup\{\epsilon\})$), $p_{0}$ the initial state,
	and $E$ the set
	of final states.

	Let $M$ be a finite automaton recognising $\mathcal{L}$.
	Suppose that $M$ has set of states $Q$,
	transition function $\delta$, initial state $q_{0}$, and set of
	final states $F$.  Extend the definition of $\delta$ so that
	$(q,\epsilon)\mapsto q$ is a valid transition for all
	$q\in Q$.

	Now define an automaton $N$ as follows.  The set of states is
	$P\times Q$, the initial state is $(p_{0},q_{0})$,
	and the set of final
	states is $E\times F$.  The transitions are defined as follows.  If
	there are transitions
	\[p_{1}\xrightarrow{d,g}p_{2}\]
	and
	\[q_{1}\xrightarrow{d}q_{2}\]
	(where $p_{1},p_{2}\in P$, $q_{1},q_{2}\in
	Q$,
	$d\in\Delta\cup\{\epsilon\}$ and $g\in\Gamma\cup\{\epsilon\}$)
	then $N$
	has a transition
	\[(p_{1},q_{1})\xrightarrow{g}(p_{2},q_{2})\]

	We prove that the new automaton recognises the set $\mathcal{L}\R$.
	Let $\beta$ be any string in $\mathcal{L}R$.  By definition of
	$\mathcal{L}\R$ we may choose
	a string $\alpha\in \mathcal{L}$ with
	$(\alpha,\beta)\in \R$.  Then we have transducer transitions
	\[p_{0}\xrightarrow{a_{1},b_{1}}p_{1}
	\xrightarrow{a_{2},b_{2}}\cdots \xrightarrow{a_{n},b_{n}}p_{n}\]
	with $p_{n}\in E$ witnessing that $(\alpha,\beta)\in R$.  Then we
	have
	$\alpha=a_{1}\ldots a_{n}$, $\beta=b_{1}\ldots b_{n}$ (where,
	possibly,
	$\epsilon$ symbols may occur).  We also have transitions
	of $M$
	\[q_{0}\xrightarrow{a_{1}}q_{2}\xrightarrow{a_{2}}\cdots
	\xrightarrow{a_{n}}q_{n}\]
	with $q_{n}\in F$ witnessing that
	$\alpha\in \mathcal{L}$.
	Then, by definition, we have transitions
	\[(p_{i-1},q_{i-1})\xrightarrow{b_{i}}(p_{i},q_{i})\]
	in $N$ demonstrating that $\beta$ is accepted by $N$.

	We reverse this argument to get the converse.  Suppose that
	$\beta\in\Gamma^{\ast}$ is accepted by $N$ via a sequence of
	transitions
	\[(p_{i-1},q_{i-1})\xrightarrow{b_{i}}(p_{i},q_{i})\]
	where $\beta=b_{1}\ldots b_{n}$ with each $b_{i}\in
	\Gamma\cup\{\epsilon\}$.  By definition of $N$ there exist
	$a_{1},\ldots,a_{n}\in \Delta\cup\{\epsilon\}$ and state transitions
	\[p_{i-1}\xrightarrow{a_{i},b_{i}}p_{i}\]
	of the transducer, and transitions
	\[q_{i-1}\xrightarrow{a_{i}} q_{i}\]
	of $M$.  This proves that $\alpha=a_{1}\ldots a_{n}\in\mathcal{L}$
	and $(\alpha,\beta)\in \R$ as required.
\end{proof}

\section{Bounded classes}\label{BoundedSection}

In this section we consider permuting machines as `black boxes' into
which input tokens are inserted and from which they eventually emerge
as output tokens.  So, at any point of a computation there may be some
tokens which are `inside' the machine (in the machine's memory)
awaiting output.  The chief
hypothesis of this section is that, for some constant $k$, the
machine can contain no more than $k$ tokens at a time (so if it is
full to capacity it must output a token before further input is
possible).  Such machines are said to be $k$-bounded.

If we consider a $k$-bounded machine as a generator of permutations
then no permutation of length $k+1$ that begins with $k+1$ can be
generated from the input $1,2,\ldots,k+1$.  Thus the closed sets
associated with $k$-bounded machines are subsets of the closed
set $\Omega_{k}$ whose basis consists of the $k!$ permutations
$k+1,a_{1},\ldots,a_{k}$ where $a_{1},\ldots,a_{k}$ ranges over all
permutations of $1,2,\ldots,k$.

We shall see shortly that permutations in $\Omega_{k}$ may be encoded as words
in a $k$-letter alphabet.  Anticipating this, we define a
subset of $\Omega_{k}$ to be regular if its encoded form is a regular
set.  We shall show that a closed subset $\mathcal{X}$
of $\Omega_{k}$ is regular
if and only if its basis
is regular.  The proof of this result is, in
principle, constructive in the sense that a recognising finite
automaton for $\mathcal{X}$ can be built from one that
recognises its basis and vice versa.  In the course of proving
this result we shall prove that it is decidable whether a regular
subset of $\Omega_{k}$ is a closed subset.

Let $\pi=\pi_{1}\pi_{2}\ldots \pi_{n}$ be a permutation of length $n$.
Its \emph{rank} encoding is the sequence
\[E(\pi)= p_1 p_2 \ldots p_n\]
where
\[p_{i} = | \{ j \mid j\geq i,\ \pi_j\leq \pi_i \} |\]
is the rank of $\pi_{i}$ among $\{\pi_{i},\pi_{i+1},\ldots,\pi_{n}\}$.

Obviously, $\pi\in\Omega_{k}$ if and only if $\pi_{1}\pi_{2}\ldots \pi_{n}$
has no
subsequence of length $k+1$ whose first element is the largest in the
subsequence and this is precisely the condition that $p_{i}\leq k$
for all $i$.  Thus every subset of $\Omega_{k}$ encodes as a subset
of $[k]^{\ast}$.

\begin{proposition}\label{regularOmega}
    $\Omega_{k}$ is regular.
\end{proposition}
\begin{proof}
    It is easy to see that a word $p=p_1p_2\ldots p_n$
    is the encoding of some permutation if and only if
    \begin{equation}
    \label{eq5}
        p_{n+1-i}\leq i \mbox{ for all } i
    \end{equation}
    (and, if this condition holds, the permutation can readily be
    calculated).
    In fact, for $p\in[k]^{\ast}$ the above inequalities may fail to hold
only for $i=1,2,\ldots, k-1$.
    Let $F$ be the set of all words of length at most $k-1$ for which
    (\ref{eq5}) does not hold.
    We now have
    \[E(\Omega_{k})=[k]^\ast \setminus [k]^\ast F,\]
    which is a regular set.
\end{proof}

\begin{example}
    Consider the closed subset $\mathcal{X}$
    of $\Omega_{2}$ whose basis is
    $312,321,231$.  The first two basis elements ensure that, indeed,
    $\X\subseteq\Omega_{2}$ so the permutations of $\mathcal{X}$
    encode as words in the alphabet $\{1,2\}$ and end with a $1$.
    It is readily checked that the third basis
    element restricts these words by prohibiting consecutive
    occurrences of the symbol $2$.

    The set of words that do contain consecutive $2$s is described by
    the regular expression $[2]^{\ast}22[2]^{\ast}$ and so is
    regular.  However $E(\mathcal{X})$ is the complement of this
    regular set within the regular set $E(\Omega_{2})$ and so is also
    regular.  Thus $\mathcal{X}$ is a regular closed set.  The
    generating function of $E(\X)$ is well-known to be
    \[\frac{1}{1-x-x^{2}}\]
    and, since $E$ is one-to-one, this is also the generating function
    of $\X$.
\end{example}

This easy example serves to illustrate that the condition of avoiding
a permutation translates into restrictions on encodings although they
are generally much more complicated than the ones above.  The argument
that proves
regularity is a very special case of more general arguments to come.

Transportation networks are another source of examples.  Theorem 1
of \cite{AtkLivTul} proves that the closed sets associated with these are
all regular.  That paper also contains an example to show that regular
closed sets need not be finitely based.

We also note that not every closed subset of $\Omega_{k}$ is regular.
Indeed, as shown in \cite{AtkMurRus},
there are uncountably many closed subsets
in $\Omega_{k}$, if $k\geq 3$; but there are
only countably many regular languages over $[k]$.

\subsection{A transducer to delete a letter from a
word}\label{deleteOne}

    Let $\pi=\pi_{1}\pi_{2}\ldots\pi_{n}$ be a permutation in
    $\Omega_{k}$ and let $p=p_{1}p_{2}\ldots p_{n}$ be its rank encoded
    form.  Let $\pi'$ be the permutation obtained from $\pi$ by
    deleting $\pi_{i}$ (and relabelling appropriately)
    and let $p'=p'_{1}\ldots p'_{i-1}p'_{i+1}\ldots
    p'_{n}$ be its encoded form.  We put
    \[\partial_{i} p = p'\]
    call this the $i$th derivative of $p$.  The process of passing from $p$ to
    $p'$ is called \emph{deleting a letter} from $p$.  We shall show how
    this may be done without ``looking at'' $\pi$.

    \begin{example}
	Let $p=2331211$ representing the permutation $\pi=2451637$.  Then
	removing the $6$th element of $\pi$ results in the permutation
	$\pi'=234156$ whose encoding is $p'=222111$.
    \end{example}

    Suppose we have to delete the $i$th letter from $p$.
    We compute $p'$ by scanning $p$ from the right.  For the positions
    to the right of $p_{i}$, each $p_{j}$ represents the rank of some
    element of $\pi$ among its successors, so these ranks will be
    unchanged by the deletion.  Therefore until we reach
    $p_{i}$ itself (which we delete) nothing changes.  But for $j<i$
    we need to know whether or not $\pi_{j}>\pi_{i}$ (so that we can
    tell
    whether or not to reduce $p_{j}$ by $1$). To do this we keep
    track of a variable $r_{j}$ defined as the rank of $\pi_{i}$ in the
    set $\{\pi_{j+1},\ldots,\pi_{n}\}$ (the number of symbols in this set
    that are less than or equal to $\pi_{i}$).  Clearly
    \[\pi_{j}>\pi_{i} \mbox{ if
	and only if }p_{j}>r_{j}\]
    Provided we have $r_{j}$ we can decide whether
    we should reduce $p_{j}$.  But, as the pointer $j$ moves to the
    left, we can easily update $r_{j}$.  Clearly, if $\pi_{j}>\pi_{i}$
    then $r_{j-1}=r_{j}$; and if $\pi_{j}<\pi_{i}$
    then $r_{j-1}=r_{j}+1$.  We therefore get Algorithm \ref{original}.

    \begin{algorithm}
    \caption{First form of the deletion algorithm}
    \begin{algorithmic}\label{original}
    \FOR{$j := n$ downto $i+1$}
        \STATE $p'_{j} := p_{j}$
    \ENDFOR
    \STATE $r_{i-1} := p_{i}$
    \FOR{$j := i-1$ downto $1$}
        \IF{$p_{j} > r_{j}$}
          \STATE $p'_{j} := p_{j}-1; r_{j-1} := r_{j}$
	  \ELSE
	  \STATE $p'_{j} := p_{j};r_{j-1} := r_{j}+1$
        \ENDIF
    \ENDFOR
    \end{algorithmic}
    \end{algorithm}

    Two easy observations make this into a finite state algorithm.
    The first is the natural programming trick to use a single
    variable $r$ in place of $r_{j}$.  The second looks odd as a
    programming trick but is nevertheless essential.  When $r_{j}\geq
    k$ the first alternative of the \textbf{if} is not followed
    nor is it followed
    thereafter; so we `freeze' $r_{j}$ to the value $k$ once it
    reaches $k$.  The result is  Algorithm \ref{modified}.

    \begin{algorithm}
    \caption{Second form of the deletion algorithm}
    \begin{algorithmic}\label{modified}
    \FOR{$j := n$ downto $i+1$}
        \STATE $p'_{j} := p_{j}$
    \ENDFOR
    \STATE $r := p_{i}$
    \FOR{$j := i-1$ downto $1$}
        \IF{$p_{j} > r$}
          \STATE $p'_{j} := p_{j}-1$
	  \ELSE
	  \STATE $p'_{j} := p_{j}$
	 \IF {$r < k$} \STATE $r := r+1$
	  \ENDIF
        \ENDIF
    \ENDFOR
    \end{algorithmic}
    \end{algorithm}

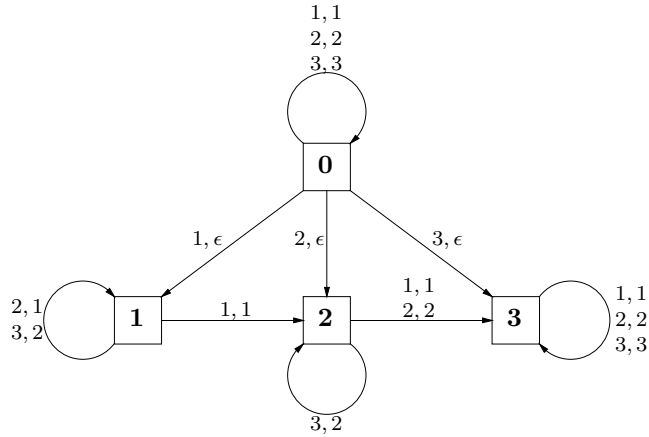
\begin{figure}[ht]
 \begin{center}
\setlength{\unitlength}{0.00041666in}
\begingroup\makeatletter\ifx\SetFigFont\undefined%
\gdef\SetFigFont#1#2#3#4#5{%
  \reset@font\fontsize{#1}{#2pt}%
  \fontfamily{#3}\fontseries{#4}\fontshape{#5}%
  \selectfont}%
\fi\endgroup%
{
\begin{picture}(8017,4950)(0,-10)
\put(3975.000,4150.000){\arc{1000.000}{2.2143}{7.2105}}
\blacken\path(4340.741,3854.777)(4275.000,3750.000)(4382.313,3811.513)(4340.741,3854.777)
\put(875.000,1500.000){\arc{1000.000}{0.6435}{5.6397}}
\blacken\path(1170.031,1865.434)(1275.000,1800.000)(1213.173,1907.133)(1170.031,1865.434)
\put(3975.000,800.000){\arc{1000.000}{5.3559}{10.3521}}
\blacken\path(3609.259,1095.223)(3675.000,1200.000)(3567.687,1138.487)(3609.259,1095.223)
\put(7075.000,1500.000){\arc{1000.000}{3.7851}{8.7813}}
\blacken\path(6779.969,1134.566)(6675.000,1200.000)(6736.827,1092.867)(6779.969,1134.566)
\path(1275,1800)(1875,1800)(1875,1200)
	(1275,1200)(1275,1800)
\path(3675,1800)(4275,1800)(4275,1200)
	(3675,1200)(3675,1800)
\path(6075,1800)(6675,1800)(6675,1200)
	(6075,1200)(6075,1800)
\path(3675,3750)(4275,3750)(4275,3150)
	(3675,3150)(3675,3750)
\path(3675,3150)(1875,1800)
\blacken\path(1953.000,1896.000)(1875.000,1800.000)(1989.000,1848.000)(1953.000,1896.000)
\path(3975,3150)(3975,1800)
\blacken\path(3945.000,1920.000)(3975.000,1800.000)(4005.000,1920.000)(3945.000,1920.000)
\path(4275,3150)(6075,1800)
\blacken\path(5961.000,1848.000)(6075.000,1800.000)(5997.000,1896.000)(5961.000,1848.000)
\path(1875,1500)(3675,1500)
\blacken\path(3555.000,1470.000)(3675.000,1500.000)(3555.000,1530.000)(3555.000,1470.000)
\path(4275,1500)(6075,1500)
\blacken\path(5955.000,1470.000)(6075.000,1500.000)(5955.000,1530.000)(5955.000,1470.000)
\put(6150,1425){\makebox(0,0)[lb]{\smash{{{\SetFigFont{10}{12}{\rmdefault}{\mddefault}{\updefault}
$\bf 3$}}}}}
\put(2350,1550){\makebox(0,0)[lb]{\smash{{{\SetFigFont{8}{9.6}{\rmdefault}{\mddefault}{\updefault}
$\begin{array}{l} 1,1 \end{array}$
}}}}}
\put(4675,1700){\makebox(0,0)[lb]{\smash{{{\SetFigFont{8}{9.6}{\rmdefault}{\mddefault}{\updefault}
$\begin{array}{l} 1,1 \\ 2, 2 \end{array}$
}}}}}
\put(-300,1425){\makebox(0,0)[lb]{\smash{{{\SetFigFont{8}{9.6}{\rmdefault}{\mddefault}{\updefault}
$\begin{array}{l} 2, 1 \\ 3, 2 \end{array}$
}}}}}
\put(1350,1425){\makebox(0,0)[lb]{\smash{{{\SetFigFont{10}{12}{\rmdefault}{\mddefault}{\updefault}
$\bf 1$
}}}}}
\put(3750,1425){\makebox(0,0)[lb]{\smash{{{\SetFigFont{10}{12}{\rmdefault}{\mddefault}{\updefault}
$\bf 2$
}}}}}
\put(3750,3375){\makebox(0,0)[lb]{\smash{{{\SetFigFont{10}{12}{\rmdefault}{\mddefault}{\updefault}
$\bf 0$
}}}}}
\put(2000,2450){\makebox(0,0)[lb]{\smash{{{\SetFigFont{8}{9.6}{\rmdefault}{\mddefault}{\updefault}
$\begin{array}{l} 1, \epsilon \end{array}$
}}}}}
\put(5050,2450){\makebox(0,0)[lb]{\smash{{{\SetFigFont{8}{9.6}{\rmdefault}{\mddefault}{\updefault}
$\begin{array}{l} 3, \epsilon \end{array}$
}}}}}
\put(3300,2450){\makebox(0,0)[lb]{\smash{{{\SetFigFont{8}{9.6}{\rmdefault}{\mddefault}{\updefault}
$\begin{array}{l} 2, \epsilon \end{array}$
}}}}}
\put(3500,100){\makebox(0,0)[lb]{\smash{{{\SetFigFont{8}{9.6}{\rmdefault}{\mddefault}{\updefault}
$\begin{array}{l} 3, 2 \end{array}$
}}}}}
\put(7375,1425){\makebox(0,0)[lb]{\smash{{{\SetFigFont{8}{9.6}{\rmdefault}{\mddefault}{\updefault}
$\begin{array}{l} 1, 1 \\ 2, 2 \\ 3, 3 \end{array}$
}}}}}
\put(3500,5000){\makebox(0,0)[lb]{\smash{{{\SetFigFont{8}{9.6}{\rmdefault}{\mddefault}{\updefault}
$\begin{array}{l} 1, 1 \\ 2, 2 \\ 3, 3 \end{array}$
}}}}}
\end{picture}
}
\caption{Deletion transducer with $k = 3$}
	 \end{center}
	     \label{diagram}

	\end{figure}

    It is now easy to define a transducer for the relation
    \[\mathcal{D}=\{(p,p')\mid p'\mbox{ is obtained by deleting one
    letter from }p\}\]
    The transducer begins in a `picking' state $0$.  Once it picks a
    letter to delete it passes through a sequence of states numbered
    according to the variable $r$ in Algorithm \ref{modified}.  The
    transducer for the case $k=3$ is shown in Figure \ref{diagram}.

    \begin{note} Strictly speaking, what we have constructed is a
    transducer for a relation where the words in question are read from
    right to left.  To avoid notational clutter we make the
    convention that all finite automata and transducers read their
    input from right to left.  Of course any conclusion that we reach
    of the form ``$\L$ is a regular language'' is independent of the
    direction of reading since $\L$ is regular if and only if its
    reverse is regular.
    \end{note}

    \begin{proposition}\label{regularSets}
	Let $\L\subseteq E(\Omega_{k})$ be regular.  Then each of the
following
	subsets is also regular, and finite automata recognising them are
	effectively computable from an automaton recognising $\L$.
	\begin{enumerate}
	    \item $\{\partial_{i}p\mid p\in\L, 1\leq i\leq |p|\}$,
	    \item $\{p\in E(\Omega_{k})\mid \partial_{i}p\in \L\mbox{, for
some }i\}$,
	    \item $\{p\in E(\Omega_{k})\mid \partial_{i}p\in\L\mbox{, for
all }i\}$.
	\end{enumerate}
    \end{proposition}
    \begin{proof}
	The first set is $\L\mathcal{D}$ and the second
	is $\L\mathcal{D}^{t}$ both of which
	are regular by Proposition \ref{technical}.  The third set is
	\begin{eqnarray*}
	    \{p\mid \partial_{i}p\not\in\L\mbox { for some }i\}^{C}&=&
	    \{p\mid \partial_{i}p\in\L^{C}\mbox { for some }i\}^{C}\\
	    &=&(\L^{C}\mathcal{D}^{t})^{C}
	\end{eqnarray*}
	Since regularity is preserved by complements the result follows again
	from Proposition \ref{technical}.
    \end{proof}

\subsection{A transducer for deleting any number of letters}

    Again let $\pi=\pi_{1}\pi_{2}\ldots\pi_{n}$ be a permutation in
    $\Omega_{k}$ and let $p=p_{1}p_{2}\ldots p_{n}$ be its rank encoded
    form.  We shall generalise the process described in the  previous
    subsection  so that it now deletes \emph{any} number
    of letters (choosing
    which ones to delete non-deterministically again).  From the
    resulting
    algorithm we shall be able to infer the existence of a transducer
    that defines the relation
    \[\mathcal{H}=\{(p,p')\mid p'\mbox{ arises by
    deleting any number of letters
    from }p\}\]
    In the generalisation a right to left scan takes place
    as before.  But now, rather than setting up a single variable $r$
    when the deleted letter is met, we have to set up a different
    variable every time we come to a letter that is to be deleted.

    So, suppose we come to a letter $p_{d}$ that we intend to
    delete.  Then we define a variable $r(d)$ (whose initial value
    will be $p_{d}$) which will play the same role as the variable $r$
    in the previous section.  Just as before when we process a letter $p_{j}$
    (either to delete it or compute the value of $p'_{j}$) we shall
    have $r(d)$ equal to the rank of $\pi_{d}$ in the
    set $\{\pi_{j+1},\ldots,\pi_{n}\}$ (that is, $r(d)$ is
    the number of symbols in this set
    that are less than or equal to $\pi_{d}$).

    Processing a particular $p_{j}$ is then done as follows:

    \begin{enumerate}
	\item if $p_{j}$ is to be deleted we set up a variable $r(j)$ as just
	mentioned and update any existing variables $r(d)$; this updating is
	explained below.
	\item if $p_{j}$ is not to be deleted we must use the variables
	$r(d)$ so far defined to compute the value of $p'_{j}$; and we must
	update these variables as necessary (see below).
    \end{enumerate}

    Exactly as before, because of the meaning of each $r(d)$ we have
    $\pi_{j}>\pi_{d}$ if and only if $p_{j}>r(d)$.  Therefore
    the number of
    $d$'s for which this occurs is the decrement that has to be applied
    to $p_{j}$ to obtain $p'_{j}$.

    To do the updating of the variable $r(d)$ (so that it has the
    appropriate value when $j$ is decreased by $1$) we notice that any $d$
    for which $p_{j}>r(d)$ means that $r(d)$ is not changed;
    otherwise it must be increased by $1$.

    The behaviour of this algorithm when a symbol $p_{j}$ is processed
    is governed by the values of the set of variables $r(d)$.  In
    order to turn the algorithm into a transducer to recognise the
    relation $\mathcal{H}$ we have to demonstrate that only a fixed
    number of variables taking a fixed set of values is required.

    First, we have the same remark as
    before: any $r(d)$ which reaches the value $k$ can never affect
    whether $p_{j}$ should be changed; so such $r(d)$'s can be discarded.  The
    second remark is that the $r(d)$ are ranks of \emph{different}
    elements within the same set ($\{\pi_{j+1},\ldots,\pi_{n}\}$);
    therefore the values $r(d)$ are distinct and
    so we shall never have more than $k-1$ of them to store.

    The state of the algorithm, as represented by the values of the
    $r(d)$, is therefore confined to one of a finite number of
    possibilities.  A convenient way of representing the state is as
    a $(0,1)$ vector $(s_{1},\ldots,s_{k-1})$.  We set
    $s_{t}=1$ if there is a variable $r(d)$ in the current `live'
    set whose value is $t$; otherwise we set $s_{t}=0$.  This
    coding of state allows the automatic `dropping'
    of a variable $r(d)$ once
    it reaches the value $k$.

    Translating the way in which the $r(d)$
    are handled, the updating of the variables $s_{t}$ when a symbol
    $p_{j}=e$ is processed is easily
    seen to be:
    \[(s_{1},\ldots,s_{k-1}):=
    (s_{1},\ldots,s_{e-1},1,s_{e},\ldots,s_{k-2})\]
    if $p_{j}$ is to be deleted and
    \[(s_{1},\ldots,s_{k-1}):=
    (s_{1},\ldots,s_{e-1},0,s_{e},\ldots,s_{k-2})\]
    otherwise.  The value of $p'_{j}$ in the latter case is
    $p_{j}-\sum_{f<e}s_{f}$.

We summarise this discussion in

\begin{proposition}
    There is a transducer that defines the relation
       \[\mathcal{H}=\{(p,p')\mid p'\mbox{ arises by
    deleting any number of letters
    from }p\}\]
\end{proposition}

\begin{figure}[ht]

        \begin{center}
\setlength{\unitlength}{0.00041666in}
\begingroup\makeatletter\ifx\SetFigFont\undefined%
\gdef\SetFigFont#1#2#3#4#5{%
  \reset@font\fontsize{#1}{#2pt}%
  \fontfamily{#3}\fontseries{#4}\fontshape{#5}%
  \selectfont}%
\fi\endgroup%
{\newcommand{\dashlinestretch}{30}
\begin{picture}(10232,7822)(0,-10)
\put(1729.500,7048.500){\arc{1497.652}{1.4765}{6.3775}}
\blacken\path(2416.992,7098.975)(2475.000,6978.000)(2536.975,7096.992)(2416.992,7098.975)
\put(8020.500,7048.500){\arc{1497.652}{3.0473}{7.9483}}
\blacken\path(7213.025,7096.992)(7275.000,6978.000)(7333.008,7098.975)(7213.025,7096.992)
\put(1729.500,757.500){\arc{1497.652}{6.1889}{11.0899}}
\blacken\path(2536.975,709.008)(2475.000,828.000)(2416.992,707.025)(2536.975,709.008)
\put(8020.500,757.500){\arc{1497.652}{4.6181}{9.5191}}
\blacken\path(7333.008,707.025)(7275.000,828.000)(7213.025,709.008)(7333.008,707.025)
\path(1875,6903)(3075,6903)(3075,5703)
	(1875,5703)(1875,6903)
\path(6675,6903)(7875,6903)(7875,5703)
	(6675,5703)(6675,6903)
\path(6675,2103)(7875,2103)(7875,903)
	(6675,903)(6675,2103)
\path(1875,2103)(3075,2103)(3075,903)
	(1875,903)(1875,2103)
\path(2475,5628)(2475,2178)
\blacken\path(2415.000,2298.000)(2475.000,2178.000)(2535.000,2298.000)(2415.000,2298.000)
\path(7275,2178)(7275,5628)
\blacken\path(7335.000,5508.000)(7275.000,5628.000)(7215.000,5508.000)(7335.000,5508.000)
\path(6600,6153)(3150,6153)
\blacken\path(3270.000,6213.000)(3150.000,6153.000)(3270.000,6093.000)(3270.000,6213.000)
\path(3150,6453)(6600,6453)
\blacken\path(6480.000,6393.000)(6600.000,6453.000)(6480.000,6513.000)(6480.000,6393.000)
\path(3150,1653)(6600,1653)
\blacken\path(6480.000,1593.000)(6600.000,1653.000)(6480.000,1713.000)(6480.000,1593.000)
\path(6600,1353)(3150,1353)
\blacken\path(3270.000,1413.000)(3150.000,1353.000)(3270.000,1293.000)(3270.000,1413.000)
\path(2925,2178)(6600,5853)
\blacken\path(6557.574,5725.721)(6600.000,5853.000)(6472.721,5810.574)(6557.574,5725.721)
\path(6824,5628)(3149,1953)
\blacken\path(3191.426,2080.279)(3149.000,1953.000)(3276.279,1995.426)(3191.426,2080.279)
\put(1425,3828){\makebox(0,0)[lb]{\smash{{{\SetFigFont{8}{9}{\rmdefault}{\mddefault}{\updefault}
$1, \lambda$}}}}}
\put(4300,1953){\makebox(0,0)[lb]{\smash{{{\SetFigFont{8}{9}{\rmdefault}{\mddefault}{\updefault}
$\begin{array}{l} 1, \lambda \\ 2, \lambda \end{array}$}}}}}
\put(2100,1375){\makebox(0,0)[lb]{\smash{{{\SetFigFont{12}{14}{\rmdefault}{\mddefault}{\updefault}
$\bf 10$}}}}}
\put(2100,6200){\makebox(0,0)[lb]{\smash{{{\SetFigFont{12}{14}{\rmdefault}{\mddefault}{\updefault}
$\bf 00$}}}}}
\put(6900,6200){\makebox(0,0)[lb]{\smash{{{\SetFigFont{12}{14}{\rmdefault}{\mddefault}{\updefault}
$\bf 01$}}}}}
\put(6900,1375){\makebox(0,0)[lb]{\smash{{{\SetFigFont{12}{14}{\rmdefault}{\mddefault}{\updefault}
$\bf 11$}}}}}
\put(4300,903){\makebox(0,0)[lb]{\smash{{{\SetFigFont{8}{9}{\rmdefault}{\mddefault}{\updefault}
$\begin{array}{l} 2, 1\end{array}$}}}}}
\put(200,1053){\makebox(0,0)[lb]{\smash{{{\SetFigFont{8}{9}{\rmdefault}{\mddefault}{\updefault}
$\begin{array}{l} 2, 1 \\ 3, \lambda \\ 3, 2 \end{array}$
}}}}}
\put(8625,1053){\makebox(0,0)[lb]{\smash{{{\SetFigFont{8}{9}{\rmdefault}{\mddefault}{\updefault}
$\begin{array}{l} 1, \lambda \\ 2, \lambda \\ 3, \lambda \\ 3, 1 \end{array}$
}}}}}
\put(200,7353){\makebox(0,0)[lb]{\smash{{{\SetFigFont{8}{9}{\rmdefault}{\mddefault}{\updefault}
$\begin{array}{l} 1, 1 \\ 2, 2 \\ 3, 3 \\ 3, \lambda \end{array}$
}}}}}
\put(4300,6753){\makebox(0,0)[lb]{\smash{{{\SetFigFont{8}{9}{\rmdefault}{\mddefault}{\updefault}
$\begin{array}{l} 2, \lambda \end{array}$
}}}}}
\put(4300,5703){\makebox(0,0)[lb]{\smash{{{\SetFigFont{8}{9}{\rmdefault}{\mddefault}{\updefault}
$\begin{array}{l} 1, 1 \\ 2, 2 \end{array}$
}}}}}
\put(8625,7353){\makebox(0,0)[lb]{\smash{{{\SetFigFont{8}{9}{\rmdefault}{\mddefault}{\updefault}
$\begin{array}{l}  2, \lambda \\ 3, \lambda \\ 3, 2 \end{array}$
}}}}}
\put(7575,3828){\makebox(0,0)[lb]{\smash{{{\SetFigFont{8}{9}{\rmdefault}{\mddefault}{\updefault}
$1, 1$
}}}}}
\put(3750,4053){\makebox(0,0)[lb]{\smash{{{\SetFigFont{8}{9}{\rmdefault}{\mddefault}{\updefault}
$1, 1$
}}}}}
\put(5175,3453){\makebox(0,0)[lb]{\smash{{{\SetFigFont{8}{9}{\rmdefault}{\mddefault}{\updefault}
$1, \lambda$}}}}}
\end{picture}
}
             \caption{Involvement transducer with $k=3$}
	     \label{InvFig}
	 \end{center}
	\end{figure}
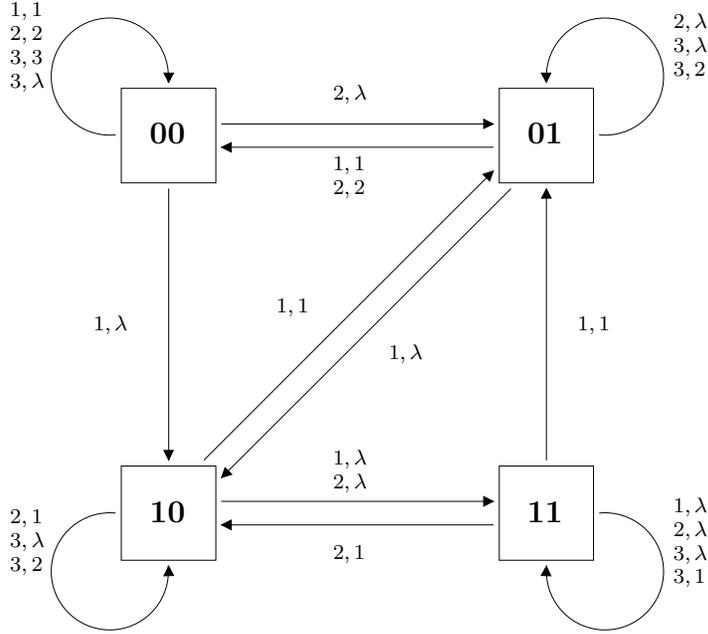

The state diagram for the transducer in the case $k=3$ is shown in
Figure \ref{InvFig}.

Clearly $\mathcal{H}^{t}$ is the relation of involvement on coded
permutations and to reflect this we write $p'\leq p$ if $p'$ can
obtained from $p$ by deleting any number of letters.

\subsection{Regularity results}

In this subsection we state and prove the main results
on $k$-bounded classes.

\begin{theorem}
\label{thm13}
There is an algorithm which decides whether or not a given regular set
$\L\subseteq [k]^\ast$ can be expressed as $\L=E(\X)$
for some closed set of permutations $\X\subseteq \Omega_k$.
\end{theorem}
\begin{proof}
First note that a set $\X$ of permutations is closed if and only if
for every $\pi=\pi_1\pi_2\ldots \pi_n\in \X$ and every $i=1,\ldots,n$,
we have $\pi\setminus \pi_i \in \X$.
Thus, $\L=E(\X)$ for some $\X$ if and only if
$\{\partial_{i}p\mid p\in\L, 1\leq i\leq |p|\}
\subseteq \L\subseteq E(\Omega_k)$.
All the three above sets are regular (Proposition \ref{regularOmega}
and Proposition \ref{regularSets}), and the automata
accepting them are known,
and hence we can decide whether these inclusions hold.
\end{proof}

\begin{theorem}
\label{RegularBasis}
A closed subset of $\Omega_{k}$ is regular if and only if its basis
is regular.
\end{theorem}

\begin{proof}
    Let $\X$ be a closed set with basis $\B$.  Suppose first that
    $\X$ is regular.  By definition
 $\B$ is the set of all permutations
$\pi=\pi_1\ldots \pi_n$ such that $\pi\not\in \X$ but
$\pi\setminus \pi_i \in \X$ for all $i=1,\ldots,n$.
Thus
$$
E(\B)=(E(\X))^C \cap \{p\mid \partial_{i}p\in E(\X)\mbox{, for all }i\},
$$
which is a regular set by Proposition \ref{regularSets}.

For the converse assume that $\B$ is regular.
    By Proposition \ref{technical} the set
    \[E(\B)\mathcal{H}^{t}=\{p\mid p'\leq
    p\mbox{ for some }p'\in E(\B)\mbox\}\]
    is regular and so its complement
    \[(E(\B)\mathcal{H}^{t})^{C}=\{p\mid p'\not\leq
    p\mbox{ for all }p'\in E(\B)\mbox\}\]
    is also regular.
    Therefore $(E(\B)\mathcal{H}^{t})^{C}\cap E(\Omega_{k})$ is
    regular as well; but this set
    is $E(\X)$ itself.
\end{proof}

The regular set operations that we have used (intersection and
complementation) are
effectively computable in the sense that automata to recognise the
resulting languages can be constructed.  Therefore we have

\begin{corollary}
There is an algorithm which, given an automaton accepting
$E(\X)$ for some regular closed set $\X$, computes
an automaton accepting $E(B)$, where $B$ is the basis of $\X$.  The
converse is also true.
\end{corollary}

This, in turn has the following pleasing consequence:

\begin{corollary}
It is decidable whether or not a given regular closed set is finitely based.
\end{corollary}

\begin{corollary}
The following are
true for any closed set $\X\subseteq \Omega_k$ with a regular
(in particular, finite) basis:
\begin{itemize}
\item[\rm (i)]
the enumeration sequence for $\X$ satisfies a linear recurrence with constant
coefficients;
\item[\rm (ii)]
membership in $\X$ can be checked in linear time.
\end{itemize}
\end{corollary}

\begin{proof}
(i) $\X$ is in one-to-one length preserving
correspondence with $E(\X)$ which, being regular, has a
rational generating function.

(ii)
Both testing for membership in a regular language and the process
of encoding permutations are linear.
\end{proof}

The first part of this corollary provides a partial (affirmative)
answer to a conjecture of Gessel (that all finitely based closed sets have
holonomic generating functions).

Theorem \ref{RegularBasis} allows us to give explicit examples of
non-regular closed sets.  Let $\mathcal{A}$ be any infinite antichain
of permutations contained in $\Omega_{k}$.  An example of such an
antichain with $k=3$ is given in \cite{AtkMurRus}.  Let
$\mathcal{A}_{0}=\{\alpha_{n_{1}},\alpha_{n_{2}},\ldots\}$ be an infinite
subset of $\mathcal{A}$ such that
\begin{enumerate}
    \item $|\alpha_{n_{i}}|=n_{i}$,
    \item $n_{1}<n_{2}<\ldots$ is not a finite union of
    arithmetic progressions.
\end{enumerate}
Then $\mathcal{A}_{0}$ is a non-regular infinite antichain and, by Theorem
\ref{RegularBasis}, defines a closed set that is not regular.

\section{Monotone segment sets}\label{MonotoneSection}

In this section we consider permuting machines with an unbounded
memory.  The memory is represented by a two-way infinite tape on which
is stored an input sequence $1,2,\ldots,n$, one token per tape square, and a
reading head moves up and down the tape.  We
consider machines $M_{\phi}$ which operate under a fixed regime of forward and
backward scans of the tape that is specified by a sequence
$\phi=f_{1}f_{2}\ldots f_{k}$  of $+$ and $-$ signs.

The machine carries out $k$ scans of the tape at the end of which all
the input symbols will have been output.  The $i$th scan is
from left to right if $f_{i}=+$ and from right to left if $f_{i}=-$.
During each scan the machine will either skip over a symbol or
output it (sequentially onto a second tape say).  Such a computation
can be described by a \emph{computation word}
$c_{1}\ldots c_{n}$ with $1\leq c_{i}\leq
k$;  the term $c_{i}$ gives the scan number on which symbol $i$ was
output.

\begin{example} \label{segmentExample}
    Let $\phi=(+,-,-)$ so that $M_{\phi}$ does one left to right
    scan and two scans right to left.  Suppose that the input tape
    contains $123456789$.  Then, supposing $M_{\phi}$
    is subject to no further constraints it might, in its first scan
    output $2,4,8$, in its second scan output $7,3$, and in its final
    scan output $9,6,5,1$.  The result is the output
    permutation $248739651$.  Notice that there is another
    computation by this machine that
    produces the same output permutation (the first scan outputs
    $2,4$, the second scan outputs $8,7,3$, and the third
    scan outputs $9,6,5,1$).  The computation words for these two
    computations are $312133213$ and $312133223$.

    Clearly this machine can only output permutations which have a
    segmentation $\alpha\beta\gamma$ where $\alpha$ is increasing
    and $\beta,\gamma$ are decreasing.  However, we do not exclude
    the possibility that, due to further constraints on the operation
    of the machine, not all permutations of this form can occur.
\end{example}

In the general case the (closed) set of permutations output by $M_{\phi}$
is a subset of
\[\W_{\phi}=\{\sigma_{1}\sigma_{2}\ldots\sigma_{k}\}\]
where each
$\sigma_{i}$ is an increasing sequence of symbols
if $f_{i}=+$ and a decreasing sequence otherwise.  The main
results of this section are that the closed subsets of
$\W_{\phi}$ have linear time recognisers and
rational generating functions.

Every computation word $c$ gives rise to a permutation $D_{\phi}(c)\in
W_{\phi}$.
To be precise, if we regard $c$ as a function
\[c:[n]\to[k]\]
then $D_{\phi}(c)$ is the permutation
obtained by concatenating the sets $c^{-1}(1)$ through $c^{-1}(k)$,
with the $i$th set in this concatenation arranged in increasing
order if $f_{i}=+$, and in decreasing order
if $f_{i}=-$.  It is easily seen that $\W_\phi$ is the image of
$[k]^{\ast}$ under the map $D_{\phi}$.

We have observed
already that $D_{\phi}$ is not one-to-one but clearly each
$D^{-1}_{\phi}(\pi)$ is a finite set (that is, every permutation
$\pi\in \W_{\phi}$
can be obtained in only finitely many ways).  We shall find it
convenient to call
its members the \emph{encodings} of $\pi$.  This situation differs
from that in the previous section in that now a permutation may have
several encodings.  Nevertheless we define subset $\X$ of $\W_{\phi}$
to be regular, if $D^{-1}_{\phi}(\X)$ is a regular
subset of $[k]^*$.

\begin{lemma}
Suppose that $s,p \in [k]^*$ and $s$ is a subword of $p$. Then
$D_{\phi}(s) \preceq D_{\phi}(p)$. Also suppose that $\sigma \preceq \pi$ are
elements of $\W_{\phi}$. Then for each encoding $p$ of $\pi$ there exists
an encoding $s$ of $\sigma$ which is a subword of $p$
\end{lemma}

\begin{proof}
The first part is immediate. For the remainder, take a subset of
the positions in $\pi$ with pattern $\sigma$. Then just take $s$
to be the subword of $p$ on the same positions.
\end{proof}

\begin{theorem}\label{regularEncoding}
Every closed subset of $\W_{\phi}$ is
regular.
\end{theorem}

\begin{proof}
Let $\X$ be a closed subset of $\W_{\phi}$ and let $B$ be its basis.
By Theorem 2.9 of \cite{AtkMurRus}
$B$ is finite.  Let $\B$ be the set of all
elements of $[k]^*$ which have a subword belonging to $D^{-1}_{\phi}(B)$.
Since $D^{-1}_{\phi}(B)$ is finite, $\B$ is regular. Suppose that $\pi
\in \X$. Then no encoding $p$ of $\pi$ can contain an element $s$
of $D^{-1}_{\phi}(B)$ as a subword, for otherwise $D_{\phi}(s) \preceq \pi$.
So $D^{-1}_{\phi}(\X) \ci \B^c$. On the other hand, if $p \in \B^c$, and
$\pi = D_{\phi}(p)$, then $\pi \in \X$ -- for if not there is some
$\sigma \in B$ with $\sigma \preceq \pi$, and then some encoding
$s$ of $\sigma$ which is a subword of $p$, a contradiction. So
$D^{-1}_{\phi}(\X) = \B^c$, which is regular.
\end{proof}

\begin{corollary}
There is a linear time recognition algorithm for any closed
subset of $\W_{\phi}$.
\end{corollary}

We cannot immediately deduce that every closed subset of
$\W_{\phi}$ has a rational generating
function since the correspondence between $\W_{\phi}$ and
$[k]^{\ast}$ is not one-to-one.  To get around this difficulty we
define, for every $\sigma\in\W_{\phi}$, a distinguished encoding
as follows.  Let $\sigma_{1}$ be the longest monotone initial
segment of $\sigma$ consistent with the sign $f_{1}$.  Having
chosen $\sigma_{1}$ we choose the next monotone segment $\sigma_{2}$
(corresponding to $f_{2}$)
also as long as possible, and we continue in this manner until all
of $\sigma$ has been segmented (necessarily with $k$ or fewer
segments).  The corresponding encoding $c_{1}\ldots c_{k}$, where
$c_{i}=j$ if $i\in\sigma_{j}$, is called the \emph{greedy} encoding of
$\sigma$.  (In the example above the first encoding was greedy, the
second was not).

\begin{lemma}
    The greedy encoding of $\W_{\phi}$ is a regular set.
\end{lemma}
\begin{proof}
    Let $p$ and $q=p+1$ be any two consecutive positions of $\phi$.
    In the greedy encoding
    $c_{1}\ldots c_{n}$ of a permutation $\sigma\in\W_{\phi}$ let
    the positions where $c_{h}=p$ be $h=i_1,i_2,\ldots,i_a$ and those
    where $c_{h}=q$ be $h=j_1,j_{2}\ldots,j_b$.  The greedy condition
    implies one of the following:
    \begin{description}
        \item[[$f_{p}=+, f_{q}=+$]]  Since $\sigma$ has adjacent segments
	$i_1,i_2,\ldots,i_a$ and $j_1,\ldots,j_b$ we have $i_a>j_1$; that is,
	in $c$,
	the final $p$ comes after the first $q$.

        \item[[$f_{p}=+, f_{q}=-$]]  Here $\sigma$ has adjacent segments
$i_1,i_2,\ldots,i_a$ and
        $j_b,\ldots,j_2,j_1$, so $i_a>j_b$; that is, the final $p$
        comes after the final $q$.

        \item[[$f_{p}=-, f_{q}=+$]]  Similarly, the first $p$ comes
        before the first $q$.

        \item[[$f_{p}=-, f_{q}=-$]]  The first $p$ comes before the
        last $q$.
    \end{description}
    Every consecutive $p,p+1$ gives a restriction on the form of a
    greedy encoding but these restrictions are all recognisable by a
    finite automaton thus completing the proof.
\end{proof}

\begin{theorem}
Every closed subset of $\W_{\phi}$ has a rational generating
function.
\end{theorem}

\begin{proof}
    Let $\X$ be any closed subset of $\W_{\phi}$.  By Theorem
    \ref{regularEncoding} $D^{-1}_{\phi}(\X)$ is regular and
    therefore $D^{-1}_{\phi}(\X)\cap\mathcal{G}$, where $\mathcal{G}$
    is the set of greedy encodings of $\W_{\phi}$, is also regular.
    But this set is in one-to-one correspondence with $\X$.
\end{proof}

\section{Final remarks}

We have shown that closed sets are the natural objects to study in the
analysis of
permuting machines.  We have also demonstrated that, when a suitable
encoding of permutations is available, finite automata are a powerful
tool in this study.  Nevertheless many problems remain.  In
particular, one natural question is how one might extend the automata
tools to use context-free encodings.  Here one might hope to prove
that certain closed sets have an algebraic generating function rather
than a rational one.  A natural candidate to
investigate are the closed subsets of $A(312)$ (which is the language
associated with a single stack) for here well-formed bracket sequences
encode permutations in a natural way.  We hope to report progress on
such problems in a subsequent paper.

Another issue is that of practicability.  The ``effective'' methods we
have developed for contructing automata frequently lead to automata
with very large numbers of states since, in particular, we often need
to convert a non-deterministic automaton to its deterministic
version.  In some special cases we have managed to contain this state
explosion and have carried out these constructions, and this gives
hope that more efficient methods may exist.

\end{document}